\documentclass[12pt]{article}

\usepackage{amsmath,amssymb}
\usepackage{graphicx}
\usepackage{psfrag}
\usepackage{epsfig}

\newtheorem{theorem}{Theorem}

\title{\bf  Colouring of Maximal $F$-free Subsets}
\author{Bill Sands}

\date{Department of Mathematics and Statistics,\\ University of Calgary,\\
Calgary, AB, Canada T2N 1N4\\ 
\bigskip\today}

\begin{document}

\maketitle

\begin{abstract}
For each finite poset $F$ with $|F|>1$, $\chi_{ac}(F)$ denotes the smallest 
integer $n$ (if it exists) such that the 
elements of every finite poset $P$ with $|P|>1$ can be coloured with at most 
$n$ colours so that every maximal $F$-free subset of $P$ with more than one 
element gets at least two colours. In this note we discuss the problem of 
determining $\chi_{ac}(F)$ for each poset $F$, give one new result, and 
summarize what is known for posets $F$ with at most four elements.
\end{abstract}

All partially ordered sets (posets) in this paper are finite. Let $F$ be a 
poset with at least two elements. A poset $P$ is $F$-{\it free} if $P$ does 
not contain a subposet 
isomorphic to $F$. A subposet $S$ of a poset $P$ is {\it maximal $F$-free} if 
$S$ is $F$-free but every subposet of $P$ strictly containing $S$ contains a 
subposet isomorphic to $F$.

\smallskip
In \cite{sands:2010}, Sands and Shen introduced the following colouring 
problem for posets. For each finite poset $F$ with $|F|>1$, let 
$\chi_{ac}(F)$ be the smallest integer $n$ (if it exists) such that the 
elements of every finite poset $P$ with $|P|>1$ can be coloured with $n$ 
colours so that every maximal $F$-free subset of $P$ with more than one 
element gets at least two colours. 

\medskip
{\it Problem}: Determine $\chi_{ac}(F)$ for each poset $F$.

\bigskip
This problem was inspired by the two special cases $F=2$ (the two-element 
antichain) and $F={\bf 2}$ (the two-element chain). A maximal 2-free subset 
of $P$ is just a maximal chain, while a maximal ${\bf 2}$-free subset of $P$ 
is a maximal antichain. In both cases, there could exist one-element maximal 
$F$-free subsets, namely any isolated point in the first case and any 
splitting point in the second, but these examples were eliminated from 
consideration in our definition of $\chi_{ac}(F)$. Thus $\chi_{ac}(2)=2$, 
attained by simply colouring all minimals of $P$ with one colour and all 
other elements of $P$ with the second colour. In contrast, it is known that 
$\chi_{ac}({\bf 2})>2$; \cite{dssw:1991} contains an example of a poset $P$ 
(with no splitting points) that cannot be two-coloured so that every maximal 
antichain gets two colours. And in fact $\chi_{ac}({\bf 2})=3$ \cite{dkt:1991}. 

\bigskip
If $F=n$ or $F={\bf n}$ for $n>2$, then $\chi_{ac}(F)$ does not exist. For 
$F=n$ we take $P=N$ (an antichain) for large $N$; then every $(n-1)$-element 
subset of $P$ is a maximal 
$n$-free subset of $P$, and at least $N/(n-2)$ colours would be needed to 
ensure that no such subset is monochromatic. Similarly for $F={\bf n}$ the 
number of colours needed is unbounded for $P$ an arbitrarily large chain. It 
is not known whether there are any other posets $F$ for which $\chi_{ac}(F)$ 
does not exist.

\bigskip
The purpose of this note is to give one new result (Theorem \ref{thm:main2}), 
and to summarize what is currently known about $\chi_{ac}(F)$ for small $F$, 
namely for $|F|\le 4$, in the hope that others can fill in some of the gaps.

\bigskip
Here is a result (Theorem 7.11) from Shen's thesis \cite{thesis:2008}:

\medskip
\begin{theorem}\label{thm:shen}
Let $F$ be a bounded poset with no interior splitting elements. Then 
$\chi_{ac}(F)\le 10$.
\end{theorem}

\bigskip
Sands and Shen \cite{sands:2010} proved

\begin{theorem}\label{thm:s&s} 
Let $F$ be a finite poset which is nonbounded (that is, does 
not have both a minimum element and a maximum element) and has no isolated 
elements. Then $\chi_{ac}(F)\le 3$.
\end{theorem}

Here is an improvement of this result for some of the above posets $F$.

\begin{theorem}\label{thm:main2}  
Suppose that 

\smallskip\noindent
(i) $F$ has at least two minimals, and 

\smallskip\noindent
(ii) every maximal of $F$ has a lower cover which is not a minimal. 

\smallskip\noindent
Then $\chi_{ac}(F)=2$.
\end{theorem}

{\it Proof}. Let $P$ be a finite poset with at least 2 elements. 
If $P$ is an antichain, then by (ii) $P$ does not contain $F$, so any 
2-colouring of 
$P$ will work. So we may suppose that $P$ contains an element $a$ of height 1. 
Then all lower covers of $a$ are minimals of $P$. Let $b$ be such a lower 
cover. Colour $P$ as follows: every element of 
$a\!\uparrow\;=\{x\in P\;|\;x\ge a\}$ is coloured 1, and all 
other elements of $P$ (including $b$) are coloured 2. We claim this works.

\smallskip
Suppose that $S$ is an $F$-free subset of $P$ coloured 1. Then $S\cup\{b\}$ is 
still $F$-free, because $b<x$ for all $x\in S$, and so by (i) no copy of $F$ in 
$S\cup\{b\}$ can contain $b$. Thus, since $b$ is coloured 2, no maximal 
$F$-free subset of $P$ can be coloured 1.

\smallskip
Suppose that $S$ is an $F$-free subset of $P$ coloured 2. Suppose that 
$S\cup\{a\}$ contained $F$. Then this copy of $F$ must contain $a$, but $a$ is 
a maximal element of $S\cup\{a\}$ (since $S$ is coloured 2), so $a$ must be 
a maximal element of $F$. By (ii), $a$ must have a lower cover which is not 
minimal in $F$, thus not minimal in $P$, which is a contradiction since $a$ 
has height 1 in $P$. Thus, since $a$ is coloured 1, no maximal 
$F$-free subset of $P$ can be coloured 2.\qquad$\Box$

\medskip
Of course, the same result is true for the duals of the posets in Theorem 
\ref{thm:main2}.

\bigskip
Of the five 3-element posets $F$, we know from \cite{sands:2010} that 
\begin{itemize}
\item $\chi_{ac}({\bf 3})=\chi_{ac}(3)=\infty$ (the chain and antichain), 
\item $\chi_{ac}(\begin{picture}(25,10)(-12,2)
\put(0,0){\line(-1,1){10}}
\put(0,0){\line(1,1){10}}
\put(0,0){\circle*{3}}
\put(10,10){\circle*{3}}
\put(-10,10){\circle*{3}}
\end{picture})=\chi_{ac}(\begin{picture}(25,10)(-3,2)
\put(10,10){\line(-1,-1){10}}
\put(10,10){\line(1,-1){10}}
\put(0,0){\circle*{3}}
\put(10,10){\circle*{3}}
\put(20,0){\circle*{3}}
\end{picture})=3$, and 
\item $\chi_{ac}(\begin{picture}(20,10)(-7,2)
\put(0,0){\line(0,1){10}}
\put(0,0){\circle*{3}}
\put(0,10){\circle*{3}}
\put(10,5){\circle*{3}}
\end{picture})\ge 4$.
\end{itemize}

Here is the current situation for the sixteen 4-element posets:

\begin{itemize}
\item $\chi_{ac}({\bf 4})=\chi_{ac}(4)=\infty$ (the chain and antichain); 
\item $\chi_{ac}(\begin{picture}(20,10)(-7,2)
\put(0,0){\line(0,1){10}}
\put(10,0){\line(0,1){10}}
\put(0,0){\circle*{3}}
\put(0,10){\circle*{3}}
\put(10,0){\circle*{3}}
\put(10,10){\circle*{3}}
\end{picture})=2$ (by \cite{sands:2010});
\item $\chi_{ac}(F)=2$ for $F=\begin{picture}(25,10)(-2,2)
\put(0,0){\line(1,1){10}}
\put(20,0){\line(-1,1){10}}
\put(10,10){\line(0,1){10}}
\put(0,0){\circle*{3}}
\put(10,10){\circle*{3}}
\put(20,0){\circle*{3}}
\put(10,20){\circle*{3}}
\end{picture}$, $\begin{picture}(25,10)(-12,2)
\put(0,0){\line(0,1){20}}
\put(-10,10){\line(1,1){10}}
\put(0,0){\circle*{3}}
\put(-10,10){\circle*{3}}
\put(0,20){\circle*{3}}
\put(0,10){\circle*{3}}
\end{picture}$, and their duals (by Theorem \ref{thm:main2});
\item $\chi_{ac}(F)\le 3$ for $F=\begin{picture}(20,10)(-7,2)
\put(0,0){\line(0,1){10}}
\put(0,0){\line(1,1){10}}
\put(10,0){\line(0,1){10}}
\put(0,0){\circle*{3}}
\put(0,10){\circle*{3}}
\put(10,0){\circle*{3}}
\put(10,10){\circle*{3}}
\end{picture}, \begin{picture}(20,10)(-7,2)
\put(0,0){\line(0,1){10}}
\put(0,0){\line(1,1){10}}
\put(10,0){\line(0,1){10}}
\put(10,0){\line(-1,1){10}}
\put(0,0){\circle*{3}}
\put(0,10){\circle*{3}}
\put(10,0){\circle*{3}}
\put(10,10){\circle*{3}}
\end{picture}, \begin{picture}(25,10)(-2,2)
\put(10,0){\line(1,1){10}}
\put(10,0){\line(-1,1){10}}
\put(10,0){\line(0,1){10}}
\put(0,10){\circle*{3}}
\put(10,10){\circle*{3}}
\put(20,10){\circle*{3}}
\put(10,0){\circle*{3}}
\end{picture}$, $\begin{picture}(25,10)(-12,2)
\put(0,0){\line(0,1){10}}
\put(-10,0){\line(1,1){10}}
\put(10,0){\line(-1,1){10}}
\put(0,0){\circle*{3}}
\put(-10,0){\circle*{3}}
\put(10,0){\circle*{3}}
\put(0,10){\circle*{3}}
\end{picture}$ (by Theorem \ref{thm:s&s});
\item $\chi_{ac}(F)\le 10$ for $F= \begin{picture}(35,10)(-18,7)
\put(0,0){\line(-1,1){10}}
\put(0,0){\line(1,1){10}}
\put(-10,10){\line(1,1){10}}
\put(10,10){\line(-1,1){10}}
\put(0,0){\circle*{3}}
\put(-10,10){\circle*{3}}
\put(10,10){\circle*{3}}
\put(0,20){\circle*{3}}
\end{picture}$ (by Theorem \ref{thm:shen});
\item  $\chi_{ac}(F)$ is unknown for $F=\begin{picture}(35,10)(-7,2)
\put(0,0){\line(0,1){10}}
\put(0,0){\circle*{3}}
\put(0,10){\circle*{3}}
\put(10,5){\circle*{3}}
\put(20,5){\circle*{3}}
\end{picture}, \begin{picture}(30,10)(-12,7)
\put(0,0){\line(0,1){20}}
\put(0,0){\circle*{3}}
\put(0,10){\circle*{3}}
\put(0,20){\circle*{3}}
\put(10,10){\circle*{3}}
\end{picture}, \begin{picture}(45,10)(-8,2)
\put(10,0){\line(1,1){10}}
\put(10,0){\line(-1,1){10}}
\put(0,10){\circle*{3}}
\put(10,0){\circle*{3}}
\put(20,10){\circle*{3}}
\put(30,0){\circle*{3}}
\end{picture},\begin{picture}(45,10)(-8,2)
\put(10,10){\line(1,-1){10}}
\put(10,10){\line(-1,-1){10}}
\put(10,10){\circle*{3}}
\put(0,0){\circle*{3}}
\put(20,0){\circle*{3}}
\put(30,0){\circle*{3}}
\end{picture}.$
\end{itemize}

\end{document}